\documentclass[10pt]{amsart} \usepackage{amsmath,amsfonts,amscd,amssymb} \oddsidemargin 1cm \evensidemargin 1cm \textwidth 14cm \textheight 20.9cm \theoremstyle{plain} \newcounter{thmcount} \newtheorem{theorem}[thmcount]{Theorem} \newtheorem{lemma}[thmcount]{Lemma} \theoremstyle{definition} \newtheorem{remark}[thmcount]{Remark} \newtheorem{example}[thmcount]{Example} \input cyracc.def \font\tencyr=wncyr10 \def\sha{\text{\tencyr\cyracc{Sh}}} \def\F{{\mathbb F}}  \def\Q{{\mathbb Q}}   \def\Z{{\mathbb Z}}  \def\A{{\mathbb A}} \def\newmathop#1{\expandafter\gdef\csname #1\endcsname{\mathop{\rm #1}\nolimits}} \newmathop{rk} \newmathop{ord} \newmathop{Gal} \newmathop{GL} \let\lar\longrightarrow \let\iso\cong \let\tensor\otimes \def\notdiv{\hbox{$\not|\,$}} \begin{document} \title{Ranks of elliptic curves in cubic extensions} \author{Tim Dokchitser} \date{October 12, 2005} \address{Tim Dokchitser\vskip 0mm Robinson College\vskip 0mm Cambridge, CB3 9AN\vskip 0mm United Kingdom} \email{t.dokchitser@dpmms.cam.ac.uk} \maketitle 
For an elliptic curve over the rationals, Goldfeld's conjecture \cite{Gol} asserts that the analytic rank $\ord_{s=1}L(E_d/\Q,s)$ of quadratic twists $E_d$ of $E$ is positive for squarefree $d$'s with density $1/2$. In other words, the analytic rank of $E$ goes up in quadratic extensions $\Q(\sqrt{d})/\Q$ half of the time. In particular, for every $E/\Q$ there are \begin{itemize} \item[(a)] infinitely many quadratic extensions where the rank goes up, and \item[(b)] infinitely many ones where it does not. \end{itemize} In fact, both (a) and (b) are known for the analytic rank and also for the algebraic rank $\rk(E/K)=\dim_{\Q}E(K)\tensor_{\Z}\Q$. 

On the other hand, root number formulas in \cite{VD,Roh} show that the situation is somewhat different for extensions $\Q(\sqrt[r]{m})/\Q$ with $r\!>\!2$ and varying $m\!>\!1$. We will be concerned with the case $r=3$, and there are examples of curves (such as $E\!=$19A3, see \cite{VD} Cor. 7) for which the analytic rank goes up in {\em every\/} non-trivial extension $\Q(\sqrt[3]{m})/\Q$; so (b) fails for cubic extensions. As for (a), the formulas do imply that the analytic rank goes up in infinitely many cubic extensions if $E/\Q$ is semistable. It turns out that the same is true of the algebraic rank for any $E$ over a number field~$K$. Thus we have 

\begin{theorem} \label{thm} Let $K$ be a number field and let $E/K$ be an elliptic curve. There are infinitely many classes $[m]\in K^*/K^{*3}\>({\sl with\ }m\in K^*)$ such that $$ \rk(E/K(\sqrt[3]{m}))>\rk(E/K)\>. $$ \end{theorem} \begin{proof} First, $E$ has finite torsion over the compositum $F=K(\mu_3,\{\sqrt[3]{m}\}_{m\in K^*})$, as every prime $v$ of $F$ has finite residue field, and prime-to-$v$ torsion injects under the reduction map modulo $v$ if $E$ has good reduction at $v$. 

Second, with $L=K(\mu_3,\sqrt[3]{m})$, the natural map $$ e: E(K)/\ell\,E(K) \lar E(L)/\ell\,E(L) $$ is injective for $\ell\ne 2,3$; in fact, the kernel-cokernel exact sequence for the Kummer maps for $E(K)$ and $E(L)$ (see \cite{Sil1} \S VIII.2) shows that $\ker e\subset H^1(\Gal(L/K),E(L)[\ell])$, which is trivial for $\ell\ne 2,3$, because the order of $\Gal(L/K)$ divides 6. 

It follows from these two facts that points of $E(K)$ can become divisible by some prime~$\ell$ only in finitely many of the extensions $K(\sqrt[3]{m})$. Thus it suffices to show that $E(L)$ is strictly larger than $E(K)$ for infinitely many distinct fields of the form $L=K(\sqrt[3]{m})$. (This argument works generally for any abelian variety and $\sqrt[r]{m}$ with $r\ge 2$). 

Now suppose $E/K$ is given by $$ E: y^2= x^3+a x+b, \qquad a,b\in K, $$ assuming for the moment that $a\ne 0$. Let $P=(x_P,y_P)$ be a non-trivial 3-torsion point on~$E$. Thus, $P$ is an inflection point, and the function~$f$ (unique up to a constant) with divisor $3(P)-3(O)$ defines a line $$ L: \> y-y_P = \kappa (x-x_P), \qquad \kappa =\frac{3x_P^2+a}{2y_P}. $$ A computation shows that $x_P=\kappa ^2/3$ and $y_P=(\kappa ^4+3a)/6\kappa$, so $L$ is defined over the field $K(\kappa )=K(x_P,y_P)=K(P)$. Parametrise $L$ by $(x,y)=(x_P\!-\!\tau/3, y_P\!-\!\kappa \tau/3)$, express the right-hand side solely in terms of $\kappa $ and $\tau$, and use this to define a map from $\A^2$ to $\A^2$. In other words, let $k$ and $t$ be indeterminants and consider the rational map $\phi: \A^2_{k,t}\to\A^2_{x,y}$ given by $$ x = \frac{k^2-t}3, \quad y=\frac{k^4+3a-2k^2t}{6k}. $$ Substituting these into the equation for $E$ shows that the Zariski closure of $\phi^{-1}(E)$ is the affine curve $$ C: \>\> 4 k^2 t^3 = k^8 + 18 a k^4 + 108 b k^2 - 27 a^2\>. $$ The degree 8 polynomial $P(x)$ on the right has discriminant $-2^{24}3^{21}a^2(4a^3+27b^2)^3\ne 0$, so $C$ is non-singular and geometrically irreducible; in fact, $C$ has geometric genus $7$. It is also clear from the construction that $P(\kappa)=0$, although the fact that the equation of $C$ has no terms with $t$ and $t^2$ is somewhat surprising, and depends on the exact choice of expressions for $x_P$ and $y_P$ in terms of $\kappa$. 

Now every $x\in K^*$ gives a point $Q_x=(x, \sqrt[3]{m_x})\in C(K(\sqrt[3]{m_x}))$ with $m_x=P(x)/4x^2$. These $Q_x$ lie in infinitely many distinct extensions $K(\sqrt[3]m)/K$, for otherwise the compositum $L=K(\{m_x\}_{x\in K^*})$ would be a number field with $C(L)$ infinite, contradicting Faltings' theorem. Finally, if $m_x\notin K^{*3}$, then the point $\phi(Q_x)$ is in $E(K(\sqrt[3]{m_x}))$ but not in $E(K)$. 

It remains to note that the same construction works when $a=0$, except that the equation of $C$ has to be divided by $k^2$, in which case $C$ has geometric genus 4 rather than~7. \end{proof} \begin{remark} For $K=\Q$, a related result due to Fearnley and Kisilevsky (\cite{FK}, Thm. 1a) is that for any $E/\Q$, the set of {\em abelian} cubic extensions $L/\Q$ for which $\rk(E/L)>\rk(E/\Q)$ is either empty or infinite. (Note also the appearance of our $P(x)$ in Prop. 3 of \cite{FK}.) \end{remark} \begin{remark} Call a prime $v$ of $K$ anomalous for $E[p]$ if $E$ has good reduction at $v$ and the reduced curve has $\tilde E(k_v)[p]\ne 0$; so $p$ is anomalous for $E/\Q$ as defined by Mazur in \cite{Maz} if it is anomalous for $E[p]$ in this terminology. 

Assume that $P(x)$ is irreducible over $K$, so that it is a minimal polynomial for~$\kappa$. Then for all but finitely many primes $v$ of $K$, $P(x)$ has a root modulo $v$ if and only if $v$ is anomalous for $E[3]$. It follows easily that apart from finitely many exceptions, every extension $K(\sqrt[3]m)/K$ produced in the proof of the theorem is ramified at some anomalous prime for $E[3]$. The appearance of anomalous primes in the construction is not coincidental, and has possibly a deep connection to Iwasawa theory. We illustrate this with one example. \end{remark} \begin{example} Take $E=X_1(11)$ of conductor 11, given by the equation $$ E: y^2=x^3-{\textstyle\frac13}x+{\textstyle\frac{19}{108}}\>. $$ If $m>1$ is a cube-free integer, then the analytic rank of $E/\Q(\sqrt[3]{m})$ is odd if and only if $11|m$. Let us look at the even rank case. 

Denote $K=\Q(\mu_3)$. From 3-descents for $E/\Q$ and $E_{-3}/\Q$, one obtains $E(K)\iso\Z/5\Z$ and $\sha(E/K)[3]=1$. It follows that the cyclotomic Euler characteristic $\chi_{cyc}(E/K)=1$, as it is the 3-part of the quantity \begin{equation} \label{chicyc} \bigl|\sha(E/K)[3^\infty]\bigr| \cdot \prod_{v|3} |\tilde E(k_v)|^2 \cdot \prod_v c_v \cdot |E(K)|^{-2}, \end{equation} and all of the terms are 3-adic units. 

Now let $F_m=K(\sqrt[3]{m})$ for some cube-free $m$, which is prime to 11. This is an abelian cubic extension, and an application of a formula by Hachimori and Matsuno for the $\lambda$-invariant in $p$-power Galois extensions shows that the following conditions are equivalent (cf. \cite{HM} Thm. 3.1 and \cite{DD} Cor. 3.20, 3.24), \begin{itemize} \item[(i)] Either $\rk(E/F_m)>0$, or $\rk(E/F_m)=0$ and $\chi_{cyc}(E/F_m)\ne 1$, \item[(ii)] $v|m$ for some prime $v$ of $K$ such that $\tilde E(k_v)[3]\ne 0$. \end{itemize} Moreover, the expression for $\chi_{cyc}(E/F_m)$ as in \eqref{chicyc} shows that (i) actually reads ``either $\rk(E/F_m)>0$ or $|\sha(E/F_m)[3]|\ne 0$'', because the other terms stay prime to 3. 

As for (ii), a prime $v$ of $K$ with $v|\ell$ ($\ell\ne3,11$) is anomalous for the 3-torsion of $E/K$ if and only if $\ell$ is anomalous for the 3-torsion of $E/\Q$. This is clear if $\ell$ splits in $K$, and for $l$ inert this follows by inspection of the possible conjugacy classes of Frobenius of $\ell\equiv2\pmod 3$ in $\Gal(\Q(E[3])/\Q)\iso\GL_2(\F_3)$. To be precise, it is not hard to see that the possible degrees of the irreducible factors of $$ P(x) = x^8 - 6x^4 + 19x^2 - 3\>, $$ modulo such $l$ are $(1,1,2,2,2)$ and $(8)$, so $P(x)$ has a linear factor over $\F_l$ if and only if it has one over $\F_{l^2}$. 

For $E=X_1(11)$, the above equivalence shows that anomalous primes are responsible for either the rank of $E$ or $\sha[3]$ going up in cubic extensions. Incidentally, this proves that $\rk(E/\Q(\sqrt[3]m))$ is zero for infinitely many $m$ (those not divisible by 11 or anomalous primes), but does not say whether it is the rank or $\sha[3]$ that goes up otherwise. On the other hand, the construction in Theorem \ref{thm} implies the following: \begin{lemma} For $E=X_1(11)$, we have $\rk(E/\Q(\sqrt[3]m))>0$ for infinitely many distinct cube-free integers $m>1$ that are prime to $11$, and infinitely many of those with $11||m$. \end{lemma} \begin{proof} It is easy to see that every $x\in\Q^*$ which is an 11-adic unit and $x\equiv\pm 1\mod 11$ (resp. $x\not\equiv\pm 1\mod 11$) gives a point $\phi(Q_x)$ of $E(\Q(\sqrt[3]{m}))$ with $11||m$ (resp. $11\notdiv m$). \end{proof} \end{example} 

\noindent {\bf Acknowledgements.} The author would like to thank J. Coates, V. Dokchitser, A. J. Scholl and J. Top for discussions of this problem.  \end{document}